\input ./style/arxiv-general.cfg
\documentclass[MSNbibl,number,citesort,nameyear,dvips]{arxbj}
\makeatletter
   \@ifpackageloaded{graphicx}{}{\usepackage{graphicx}}
\makeatother
\usepackage{dcolumn}

\aid{0}
\volume{23}
\issue{1}
\pubyear{2017}
\firstpage{1}
\lastpage{22}
\doi{10.3150/14-BEJ671} 
\docsubty{FLA}

\makeatletter

\newcommand{\lleft}{\left}
\newcommand{\rright}{\right}

\newtheorem{theorem}{Theorem}[section]
\newtheorem{lemma}{Lemma}[section]
\newtheorem{proposition}{Proposition}[section]
\newtheorem{corollary}{Corollary}[section]
\newremark{ChromosomeRepresentation}{Chromosome representation}
\newremark{IslandModelandConvergence}{Island Model and convergence}
\newcommand{\ignore}[1]{}
\newcommand{\newtext}[1]{{#1}}

\newcommand{\eqref}[1]{(\ref{#1})}

\newcommand{\QAR}{\operatorname{QAR}}
\newcommand{\AR}{\operatorname{AR}}

\newcolumntype{d}[1]{D{.}{.}{#1}}
\makeatother

\begin{document}
\begin{frontmatter}

\title{Piecewise quantile autoregressive modeling for nonstationary
time series}
\runtitle{Piecewise quantile autoregressive modeling}

\begin{aug}
\author[A]{\inits{A.}\fnms{Alexander}~\snm{Aue}\corref{}\thanksref{e1}\ead[label=e1,mark]{aaue@ucdavis.edu}},
\author[A]{\inits{R.C.Y.}\fnms{Rex~C.Y.}~\snm{Cheung}\thanksref{e2}\ead[label=e2,mark]{rccheung@ucdavis.edu}},
\author[A]{\inits{T.C.M.}\fnms{Thomas~C.M.}~\snm{Lee}\thanksref{e3}\ead[label=e3,mark]{tcmlee@ucdavis.edu}}
\and
\author[A]{\inits{M.}\fnms{Ming}~\snm{Zhong}\thanksref{e4}\ead[label=e4,mark]{mgzhong@ucdavis.edu}}
\address[A]{Department Statistics,
University of California,
One Shields Avenue,
Davis, CA 95616, USA.\\ \printead{e1,e2,e3,e4}}
\end{aug}

\received{\smonth{1} \syear{2014}}
\revised{\smonth{7} \syear{2014}}

%
\begin{abstract}
We develop a new methodology for the fitting of nonstationary time
series that exhibit nonlinearity, asymmetry, local persistence and
changes in location scale and shape of the underlying distribution. In
order to achieve this goal, we perform model selection in the class of
piecewise stationary quantile autoregressive processes. The best model
is defined in terms of minimizing a minimum description length
criterion derived from an asymmetric Laplace likelihood. Its practical
minimization is done with the use of genetic algorithms. If the data
generating process follows indeed a piecewise quantile autoregression
structure, we show that our method is consistent for estimating the
break points and the autoregressive parameters. Empirical work suggests
that the proposed method performs well in finite samples.
\end{abstract}

%
\end{frontmatter}

\section{Introduction}

Many time series observed in practice display nonstationary behavior,
especially if data is collected over long time spans. Nonstationarity
can affect the trend, the variance--covariance structure or, more
comprehensively, aspects of the underlying distribution. Since
estimates and forecasts can be severely biased if nonstationarity is
not properly taken into account, identifying and locating structural
breaks has become an important issue in the analysis of time series.
Over the years, there has been a large amount of research on issues
related to testing and estimating structural breaks in sequences of
independent random variables, time series and regression models. Most
of these focus on considering breaks in the (conditional) mean, while a
smaller number of publications are available for breaks in the
(conditional) variance. The relevant lines of research are summarized
in the monograph \cite{Csorgo:Horvath:1997} and the more recent
survey paper \cite{Aue:Horvath:2013}.

In various situations, however, it may be helpful and more informative
to study structural breaks in the (conditional) quantiles. As a case in
point, Hughes \textit{et al.} \cite
{Hughes:SubbaRao:SubbaRao:2007} have argued convincingly
that the increase in mean surface temperatures recorded at temperature
stations across the Antarctic can to a large degree be attributed to an
increase in the minimum and lower quantile temperatures. When focusing
on the mean, this additional information about the underlying changes
in variation is smoothed out and unavailable for a more in-depth
analysis. As another example, the Value at Risk, a measure of loss
associated with a rare event under normal market conditions, is by
definition a quantile and more important for risk managers than
information on measures of central tendency such as the mean.

Global estimation procedures for quantiles are often performed in the
quantile regression framework described in \cite{Koenker:2005}. There
is by now a rich body of literature on the various aspects of quantile
regression models. Detecting structural breaks in nonstationary time
series over different quantiles, however, is a comparatively new
research area. Contributions in a different direction from ours include
\cite{Bai:1998}, who considered the estimation of structural breaks
in the median of an underlying regression model by means of least
absolute deviations. In the quantile regression framework, Aue \textit{et~al.} \cite{Aue:Cheung:Lee:Zhong:2014} have
recently developed a related
methodology to perform segmented variable selection that includes break
point detection as a special case. The focus of the present paper,
however, is more on the aspects of nonlinear time series analysis.

In order to capture nonlinearities such as asymmetries, local
persistence, and changes in location, scale and shape, in conjunction
with temporal dependence that is frequently observed in applications,
and thus to obtain a more complete picture of the distributional
evolution of the underlying random processes, we propose in this paper
a new method for estimating structural breaks at any single quantile or
across multiple quantiles. Our methodology differs from the works above
in that it is not based on hypothesis testing. Instead we try to match
the observed data with a best fitting piecewise quantile
autoregression. These models, introduced by Koenker and Xiao \cite{Koenker:Xiao:2006},
are members of the class of random coefficient autoregressions that
allow the autoregressive coefficients to be quantile dependent and,
therefore, generalize linear quantile autoregressions as studied by
Koul and Saleh \cite{Koul:Saleh:1995}, and Hallin and Jure{\v{c}}kov{\'a} \cite{Hallin:Jureckova:1989}, among
others. We discuss quantile autoregression models and their piecewise
specifications in Section~\ref{sec:model}. In particular, we state
necessary and sufficient conditions for the existence of stationary
solutions and discuss the estimation of the parameters via optimizing a
subgradient condition. These results will then be generalized to the
piecewise stationary case.

Recognizing the connection between estimation of quantile
autoregression parameters and maximum likelihood estimation for
asymmetric Laplace random variables \cite{Yu:Lu:Stander:2003}, we
shall apply the minimum description length principle \cite
{Rissanen:1989} to define the best fitting piecewise quantile
autoregression. Details of this are given in Section~\ref{sec:mdl}.
Minimization of the resulting convex objective function will then yield
the best fitting model for the given data. The numerical complexity of
this optimization problem is handled via the application of a genetic
algorithm \cite{Davis:1991}.

From a technical perspective, our methodology is related to \cite
{Davis:Lee:Rodriguez-Yam:2006}, who proposed an automatic procedure
termed Auto-PARM. This procedure is designed to detect structural
breaks by fitting piecewise stationary, linear autoregressive time
series models which are estimated through the minimization of a minimum
description length criterion using a normal likelihood. Auto-PARM is
defined to mimic the second-order properties of the data but is not
always able to adjust to a nonlinear framework and does not provide
additional insight into distributional changes other than those
affecting the conditional mean and variance of the data given past observations.

The remainder of the paper is organized as follows. In Section~\ref
{sec:model}, quantile autoregressive models are introduced. Estimation
and model selection aspects for piecewise quantile autoregressive
models are detailed in Section~\ref{sec:mdl}. Sections~\ref
{sec:theory} and \ref{sec:ga} deal with asymptotic results and
implementation details, respectively. Empirical properties of the
proposed methodology are evaluated through simulations in Section~\ref
{sec:sim} and
real data examples in Section~\ref{sec:app}. Section~\ref{sec:sum}
concludes and
all technical proofs are given in the \hyperref[sec:proof]{Appendix}.

\section{Quantile autoregressions}
\label{sec:model}


Linear autoregressive models have played a dominant role in classical
time series analysis for at least half a century. The popularity stems
partially from their closeness to the linear regression framework with
its well-developed theory. They are, however, unable to capture
nonlinear dynamics and local persistence. With the objective of
dynamically modeling the evolution of location, scale and shape of the
underlying processes, Koenker and Xiao \cite
{Koenker:Xiao:2006} have introduced a
particular subclass of random coefficient autoregressive models called quantile
autoregressions. In this model, autoregressive coefficients are allowed
to vary with the quantiles $\tau\in[0,1]$. In contrast to many of the
standard contributions to the random coefficient autoregression area
for which independence is a key assumption, the coefficients possess a
strong functional
relationship; in sequel $\mathbb{Z}$ denotes the set of integers. A
real time series $(y_t\dvt t\in\mathbb{Z})$ is said to follow a
quantile autoregression of order $p$, shortly $\QAR(p)$, if
%
%
\begin{equation}
\label{eqn:qar1} y_{t}=\theta_{0}(u_t)+
\theta_{1}(u_t)y_{t-1}+\cdots+\theta
_{p}(u_t)y_{t-p}, \qquad t\in\mathbb{Z},
\end{equation}
where $(u_t\dvt t\in\mathbb{Z})$ are independent random variables
distributed uniformly on the interval $[0,1]$, and $\theta_{j}\dvtx
[0,1]\to\mathbb{R}$, $j=0,1,\ldots,p$, are the coefficient
functions. In order to exhibit the connection to standard random
coefficient autoregressions, (\ref{eqn:qar1}) can also be written more
conventionally in the form
%
%
\begin{equation}
\label{eqn:qar4} y_{t}=\phi_{0}+\phi_{1,t}y_{t-1}+
\cdots+\phi_{p,t}y_{t-p}+\varepsilon_{t},\qquad t\in
\mathbb{Z},
\end{equation}
where $\phi_{0}=E\{\theta_{0}(u_t)\}$, $\varepsilon_{t}=\theta
_{0}(u_t)-\phi_{0}$, and $\phi_{j,t}=\theta_{j}(u_t)$ for
$j=1,\ldots,p$ and $t\in\mathbb{Z}$. We have in particular that the
innovations $(\varepsilon_{t}\dvt t\in\mathbb{Z})$ constitute an
independent, identically distributed sequence with distribution
function $F(\cdot)=\theta_{0}^{-1}(\cdot+\phi_{0})$. Therefore,
necessary and sufficient conditions for the existence of a strictly
stationary solution to the equations (\ref{eqn:qar1}) can be derived
from the work of Aue \textit{et al}. \cite
{Aue:Horvath:Steinebach:2006}, which also
contains statements concerning the finiteness of moments of quantile
autoregressions.


The estimation of the quantile autoregression functions $\theta(\tau
)$ in stationary quantile autoregressive models (\ref{eqn:qar1}) is
typically achieved \cite{Koenker:2005} by solving the convex
optimization problem
%
%
\begin{equation}
\label{eqn:qar3} \min_{\theta(\tau)\in\mathbb{R}^{p+1}}\sum_{t=1}^{n}
\rho_{\tau
}\bigl\{y_{t}-X'_{t}\theta(
\tau)\bigr\},
\end{equation}
where $\rho_{\tau}(u)=u\{\tau-I(u<0)\}$ is the check function.
Solutions $\hat\theta(\tau)$ of (\ref{eqn:qar3}) are called
autoregression quantiles. Asymptotic properties of the estimation
procedure have been derived in \cite{Koenker:Xiao:2006}. It should be
noted that the assumptions under which the following proposition holds
require $X_t^\prime\theta(\tau)$ to be monotonic.
This will not always be reasonable. However, for the methodology
developed in this paper, this is not an issue insofar as we derive
asymptotic statements only about the quality of the segmentation
procedure but not on the quality of the estimator $\hat\theta$.
%
%
\begin{proposition}\label{prop:3}
Let $F_{t-1}={P}(y_t<\cdot\mid\mathcal{F}_{t-1})$ be the conditional
distribution function of $y_t$ given $\mathcal{F}_{t-1}$, and denote
by $f_{t-1}$ its derivative. Under stationarity and if $f_{t-1}$ is
uniformly integrable on $\mathcal{X}=\{x\dvt0<F(x)<1\}$, then
\[
\Sigma^{-1/2}n^{1/2}\bigl[\hat\theta(\cdot)-\theta(\cdot)\bigr]
\stackrel{\mathcal{D}} {\longrightarrow}B_{p+1}(\cdot)\qquad (n\to
\infty),
\]
where $\Sigma=\Omega_1^{-1}\Omega_0\Omega_1^{-1}$ with $\Omega
_0=E(X_tX_t^\prime)$ and $\Omega_1=\lim_n\frac{1}n\sum
_{t=1}^nf_{t-1}\{F_{t-1}^{-1}(\tau)\}X_tX_t^\prime$. Moreover,
$(B_{p+1}(\tau)\dvt\tau\in[0,1])$ is a standard
$(p+1)$-dimensional Brownian bridge.
\end{proposition}

If the number of break points $m$ is given, then estimating their
locations and the $m+1$ piecewise quantile autoregressive models at a
specific quantile $\tau\in(0,1)$ can be done via solving
%
%
\begin{equation}
\label{eqn:qar7} \min_{\theta(\tau),\mathcal{K}}\sum_{j=1}^{m+1}
\sum_{t=k_{j-1}+1}^{k_{j}}\rho_{\tau}\bigl
\{y_{t}-X'_{j,t}\theta_{j}(\tau)
\bigr\}.
\end{equation}
Given that the number of observations in each segment increases as a
fraction of the overall sample size, the limit behavior of (\ref
{eqn:qar7}) follows directly from Proposition~\ref{prop:3}. \newtext
{For unknown $m$, we use a model selection approach to select the
numbers of segments.} To this end, we discuss the relation between
(\ref{eqn:qar3}) and (\ref{eqn:qar7}), and optimizing the likelihood
obtained from asymmetric Laplace distributions next.

The connection between the asymmetric Laplace distribution and quantile
regression has long been recognized and has often been used in the
Bayesian context. Yu \textit{et al.} \cite
{Yu:Lu:Stander:2003} have made this explicit.
If we assume that at the $\tau$th quantile the innovations
$(\varepsilon_t\dvt t\in\mathbb{Z})$ in model (\ref{eqn:qar4})
follow an asymmetric Laplace distribution with parameter $\tau$, then
maximizing the likelihood function
\[
L\bigl\{\theta(\tau)\bigr\}\propto\exp\Biggl[-\sum
_{t=1}^n\rho_\tau\bigl\{
y_t-X_t^\prime\theta(\tau)\bigr\} \Biggr]
\]
is equivalent to solving the problem in (\ref{eqn:qar3}). The
equivalent to (\ref{eqn:qar7}) could be stated in a similar fashion.
The use of the asymmetric Laplace likelihood allows us to formulate a
minimum description length criterion in order to do model selection
with (\ref{eqn:qar7}).

\section{Piecewise quantile autoregressions}
\label{sec:mdl}

\subsection{The model}
\label{subsec:piecewise qar}

Koenker and Xiao \cite{Koenker:Xiao:2006} have pointed
out that a fitted quantile
autoregressive model should serve as a useful local approximation to a
potentially more complicated global dynamic. While a single quantile
autoregression fit can already adequately and quite explicitly describe
local persistence and seemingly explosive behavior (see Sections~\ref
{sec:sim} and \ref{sec:app} for examples), it does not provide us with
means to fit nonstationary data. We propose to match a nonstationary
time series by blocks of different stationary quantile autoregressions.

The piecewise stationary quantile autoregressive models are defined as
follows. Assume that the data $y_1,\ldots,y_n$ can be segmented into
$m+1$ stationary pieces, and that, for $\ell=1,\ldots,m+1$, the $\ell
$th piece can be modeled by a $\QAR(p_\ell)$ process. For $\ell
=1,\ldots,m+1$, we denote by $k_\ell$ the $\ell$th break date, that
is, the time lag at which the transition from the $\ell$th to the
$(\ell+1)$th segment occurs. Using the convention $k_0=1$ and
$k_{m+1}=n$ and letting $u_1,\ldots,u_n$ be independent standard
uniform random variables, the $\ell$th segment is, for $t=k_{\ell
-1}+1,\ldots,k_\ell$, given by
%
%
\begin{equation}
\label{eqn:qar5} y_{t}=\theta_{\ell,0}(u_t)+
\theta_{\ell,1}(u_t)y_{t-1}+\cdots+
\theta_{\ell,p_\ell}(u_t)y_{t-p_\ell} =X_{\ell,t}^\prime
\theta_\ell(u_t),
\end{equation}
where $X_{\ell,t}=(1,y_{t-1},\ldots,y_{t-p_{\ell}})^\prime$ and
$\theta_{\ell}(u_t)=\{\theta_{\ell,0}(u_t),\ldots,\theta_{\ell
,p_{\ell}}(u_t)\}^\prime$. At $\tau\in(0,1)$, model (\ref
{eqn:qar5}) is determined by the parameters $m$, $\mathcal
{K}=(k_1,\ldots,k_m)^\prime$ and $\theta(\tau)=\{\theta_1(\tau
)^\prime,\ldots,\break \theta_{m+1}(\tau)^\prime\}^\prime$, where the
segment autoregression functions are denoted by $\theta_\ell(\tau)=\{
\theta_{\ell,0}(\tau),\theta_{\ell,1}(\tau),\ldots,\allowbreak \theta_{\ell
,p_\ell}(\tau)\}^\prime$. Observe that in the case that $m=0$, (\ref
{eqn:qar5}) reduces to the single $\QAR(p)$ model (\ref{eqn:qar1}). One
can fit the model (\ref{eqn:qar5}) even if it is not the true data
generating process and that we can then view the piecewise quantile
autoregressive structure as an approximation.

The approach taken in this paper is related to the piecewise AR model
fitting technique Auto-PARM developed in \cite
{Davis:Lee:Rodriguez-Yam:2006}. These authors utilized linear time
series models, changing the coefficient functions $\theta_{\ell
,j}(\cdot)$ in \eqref{eqn:qar5} to constants, say, $\phi_{\ell,j}$,
and were concerned mainly about matching the second-order structure of
the data with stationary AR segments. The present paper focuses on
nonlinear aspects of the time series as observed from quantiles,
thereby enabling a more comprehensive study of changes in the
distribution of the underlying data. The switch from linear to
nonlinear time series means in particular that somewhat different
arguments are needed in order to prove large-sample results (see
Section~\ref{sec:theory}). In terms of practical estimation, the
genetic algorithm behind Auto-PARM can be modified for the piecewise
quantile autoregression fitting. Details are given in Section~\ref{sec:ga}.

\subsection{Model selection at a single quantile}

In this section, we derive a minimum description length criterion for
choosing the best fitting model from the piecewise quantile
autoregressive models defined in (\ref{eqn:qar5}). As to be seen
below, the ``best'' model is defined as the one that enables the best
compression of the observed series $Y=(y_1,\ldots,y_n)^\prime$. For
introductory material on this, see, for example, \cite
{Rissanen:1989,Hansen:Yu:2000,Lee:2001}.

There are different versions of the minimum description length
principle, and the version adopted here is the so-called two-part code.
It begins with splitting $Y$ into two parts. The first part, denoted by
$\hat{\mathcal{{F}}}$, represents the fitted piecewise quantile
autoregression, and the second\vspace*{2pt} part, denoted by $\hat{\mathcal
{E}}=Y-\hat{Y}$, represents the residuals, where $\hat{Y}$ is the
fitted value for $Y$. Notice that once $\hat{\mathcal{{F}}}$ and
$\hat{\mathcal{E}}$ are known, $Y$ can be completely retrieved. The
idea of the minimum\vspace*{2pt} description length principle is to find the best
pair of $\hat{\mathcal{{F}}}$ and $\hat{\mathcal{E}}$ so that via
encoding (or compressing) $\hat{\mathcal{{F}}}$ and $\hat{\mathcal
{E}}$, $Y$ can be transmitted (or stored) with the least amount of
codelength (or memory). To quantify this idea, let $\textsc
{cl}_{\mathcal{F}}(Z|\tau)$ denote the codelength of an object $Z$
using model $\mathcal{F}$ at a specific quantile~$\tau$. Then we have
the decomposition
%
%
\begin{equation}
\label{eqn:twopart} \textsc{cl}_{\mathcal{F}}(Y|\tau)=\textsc
{cl}_{\mathcal{F}}(
\hat{\mathcal{{F}}}|\tau)+\textsc{cl}_{\mathcal{F}}(\hat{\mathcal
{E}}|\hat{
\mathcal{{F}}},\tau)
\end{equation}
for the data $Y$. In the above $\textsc{cl}_{\mathcal{F}}(Y|\tau)$
is the codelength for $Y$, $\textsc{cl}_{\mathcal{F}}(\hat{\mathcal
{{F}}}|\tau)$ is the codelength for $\hat{\mathcal{{F}}}$, while
$\textsc{cl}_{\mathcal{F}}(\hat{\mathcal{E}}|\hat{\mathcal
{{F}}},\tau)$ is the codelength for $\hat{\mathcal{E}}$. The minimum
description length\vspace*{1pt} principle defines the best fitting $\hat{\mathcal
{F}}$ as the one that minimizes $\textsc{cl}_{\mathcal{F}}(Y|\tau)$.

Using the estimated quantile autoregression structure, we obtain the
following expression:
%
%
\begin{eqnarray}
\label{eqn:qar8} \textsc{cl}_{\mathcal{F}}(\hat{\mathcal
{{F}}}|\tau)
&=&\textsc{cl}_{\mathcal{F}}(m|\tau)+\textsc{cl}_{\mathcal
{F}}(k_{1},
\ldots,k_{m}|\tau)+\textsc{cl}_{\mathcal
{F}}(p_{1},
\ldots,p_{m+1}|\tau) \nonumber\\
&&{}+\textsc{cl}_{\mathcal{F}}\bigl\{\hat{
\theta}_{1}(\tau),\ldots,\hat\theta_{m+1}(\tau)\bigr\}
\nonumber
\\[-8pt]\\[-8pt]
&=&\textsc{cl}_{\mathcal{F}}(m|\tau)+\textsc{cl}_{\mathcal
{F}}(n_{1},
\ldots,n_{m+1}|\tau)+\textsc{cl}_{\mathcal
{F}}(p_{1},
\ldots,p_{m+1}|\tau)\nonumber\\
&&{} +\textsc{cl}_{\mathcal{F}}\bigl\{\hat{
\theta}_{1}(\tau),\ldots,\hat\theta_{m+1}(\tau)\bigr\}.
\nonumber
\end{eqnarray}
To proceed further, we need the following coding result: the codelength
for an integer $T$ is $\log_2 T$ bits, leading to $\textsc
{cl}_{\mathcal{F}}(m|\tau)=\log_2 m$ and $\textsc{cl}_{\mathcal
{F}}(p_{1},\ldots,p_{m+1}|\tau)=\sum_{j=1}^{m+1}\log_{2}p_{j}$. On
the other hand, if the upper bound $T_U$ of an integer $T$ is known,
the corresponding codelength is $\log_2 T_U$ bits. This gives $\textsc
{cl}_{\mathcal{F}}(n_{1},\ldots,n_{m+1}|\tau)=(m+1)\log_2 n$, as
each $n_\ell$ is upper-bounded by $n$. Lastly, Rissanen \cite{Rissanen:1989}
has shown that \newtext{a} maximum likelihood estimate computed from
$n$ data points can be effectively encoded with $\frac{1}{2}\log_2 n$
bits. Applying this to the $\hat\theta_{\ell}(\tau)$'s, we have
$\textsc{cl}_{\mathcal{F}}\{\hat{\theta}_{1}(\tau),\ldots,\hat
\theta_{m+1}(\tau)\}=\newtext{\sum_{j=1}^{m+1}\frac
{p_{j}+1}{2}\log_{2}n_{j}}$. Combining these codelength expressions,
(\ref{eqn:qar8}) becomes
%
%
\begin{equation}
\label{eqn:qar8a} \textsc{cl}_{\mathcal{F}}(\hat{\mathcal
{{F}}}|\tau) =
\log_{2}m+(m+1)\log_{2}n+\sum_{j=1}^{m+1}
\log_{2}p_{j}+\newtext{\sum_{j=1}^{m+1}
\frac{p_{j}+1}{2}\log_{2}n_{j}}.
\end{equation}

Now for the last term in (\ref{eqn:twopart}). It is shown in \cite
{Rissanen:1989} that the codelength of the residuals $\hat{\mathcal
{E}}$ is the negative of the log likelihood of the fitted model $\hat
{\mathcal{{F}}}$. Utilizing the asymmetric Laplace likelihood this
leads to
%
%
\begin{equation}
\label{eqn:qar9} \textsc{cl}_{\mathcal{F}}(\hat{\mathcal{E}}|\hat
{\mathcal
{{F}}},\tau)= -\log L\bigl\{\theta(\tau)\bigr\} =\sum
_{j=1}^{m+1}\sum_{t=k_{j-1}+1}^{k_{j}}
\rho_\tau(\newtext{\hat{\varepsilon}_t})-n\log\bigl\{
\tau(1-\tau)\bigr\}.
\end{equation}
Combining equations (\ref{eqn:twopart}), (\ref{eqn:qar8}) and (\ref
{eqn:qar9}) and dropping the constant term $-n\log\{\tau(1-\tau)\}$,
we define the best fitting piecewise quantile autoregressive model at a
single quantile $\tau\in(0, 1)$ as the one that minimizes the minimum
description length criterion
%
%
\begin{eqnarray}
\label{eqn:qar10}
&&\textsc{mdl}(m,k_{1},\ldots,k_{m},p_{1},
\ldots,p_{m+1}|\tau)\nonumber\\
&&\quad =\log_{2}m+(m+1)\log_{2}n\\
&&\qquad{}+\sum
_{j=1}^{m+1}\log_{2}p_{j}+
\newtext{\sum_{j=1}^{m+1}\frac{p_{j}+1}{2}
\log_{2}n_{j}} 
+\sum
_{j=1}^{m+1}\sum_{t=k_{j-1}+1}^{k_{j}}
\rho_\tau(\newtext{\hat{\varepsilon}_{t}}).
\nonumber
\end{eqnarray}

\subsection{Model selection at multiple quantiles}

To extend the scope of detecting break points at a single quantile, it
is worthwhile to study the joint estimation of, say, $L$ quantiles in
order to gain more insight into the global behavior of the process. To
estimate break points for multiple quantiles, it can, for example, be
assumed that the true break locations are the same across the different
quantiles under consideration. This could lead to a borrowing of
strength in the segmentation procedure because information on the
behavior of various quantiles is added into the analysis. Instead of
summing up the minimum description length function defined in (\ref
{eqn:qar10}) for all $L$ quantiles, one could also use their weighted
sums. That is,
%
%
\begin{eqnarray}
\label{eqn:qar11} &&\textsc{mdl}(m,k_{1},\ldots,k_{m},p_{1},
\ldots,p_{m+1}|\tau_{1},\ldots,\tau_{L})
\nonumber\\[-8pt]\\[-8pt]
&&\quad =\sum_{\ell=1}^{L}\omega_{\ell}
\textsc{mdl}(m,k_{1},\ldots,k_{m},p_{1},
\ldots,p_{m+1}|\tau_{\ell}).
\nonumber
\end{eqnarray}
The weights can either be chosen in advance or data-adaptively. In the
latter case it may be worthwhile to read the discussion in Chapter~5.5
of \cite{Koenker:2005}, where similar ideas are discussed in a
location-shift regression model. For this case the optimal weights
$\omega_\mathrm{opt}=(\omega_{1,\mathrm{opt}},\ldots,\omega
_{L,\mathrm{opt}})^\prime$ are given by
$
\omega_{\mathrm{opt}}=W^{-1}v$,
where $W$ is the $L\times L$ matrix with entries $A_{\ell,\ell^\prime
}=\min\{\tau_\ell,\tau_{\ell^\prime}\}-\tau_\ell\tau_{\ell
^\prime}$ and $v=(v_1,\ldots,v_L)^\prime$ with $v_\ell
=f(F^{-1}(\tau_\ell))$. For the more complicated model under
consideration here, one could use these results as a starting point for
a more detailed analysis.

On the other hand, one could also think about a more general version of
the segmentation procedure that would not enforce simultaneous breaks
across the quantiles under consideration. Such an approach may be
useful if it could be coupled with prior information on the effect
breaks would have on the underlying distribution; for example, if
breaks would propagate in a monotone way from the lower to the upper
quantiles. The resulting minimum description length criterion would
then be even more complex. While a few issues concerning multiple
quantiles are highlighted in the empirical parts of the paper, any
detailed analysis of such modeling is, however, beyond the scope of the
present paper.

\section{Large sample results}
\label{sec:theory}

To study large sample properties assume that the underlying true model
indeed follows the piecewise quantile autoregressive structure in (\ref
{eqn:qar5}). We denote the true number of break points and their
locations respectively by $m^0$ and $k_{j}^{0}$, $j=1,\ldots,m^0$,
where $k_{j}^{0}=\lfloor\lambda^{0}_{j}n\rfloor$ and $0<\lambda
^{0}_{1}<\lambda^{0}_{2}<\cdots<\lambda^{0}_{m^0}<1$. Following
standard convention in order to ensure sufficient separation of the
break points, we choose an $\epsilon>0$ such that $\epsilon\ll\min
_{j=1,\ldots,m^0+1}(\lambda_{j}^{0}-\lambda_{j-1}^{0})$ and set
\[
\Lambda_{m}= \bigl\{(\lambda_{1},\ldots,
\lambda_{m})\dvt0<\lambda_{1}<\cdots<
\lambda_{m}<1, \lambda_{j}-\lambda_{j-1}\geq
\epsilon, j=1,2,\ldots,m+1 \bigr\},
\]
where $\lambda_{0}=0$ and $\lambda_{m+1}=1$. Fix $\tau\in(0,1)$,
and set $\lambda=(\lambda_{1},\ldots,\lambda_{m})$ and
$p=(p_{1},\ldots,p_{m+1})$. The parameters $m$, $\lambda$ and $p$ are
estimated by minimizing the minimum description length criterion
%
%
\begin{equation}
\label{mdl-lim} (\hat{m},\hat\lambda,\hat{p})=\arg\min
_{(m,\lambda,p)\in
\mathcal{M}}
\frac{1}{n}\textsc{mdl}(m,\lambda,p|\tau),
\end{equation}
where the minimum is taken in the set $\mathcal{M}=\{(m,\lambda
,p)\dvt m\leq M_0, \lambda\in\Lambda_m, 0\leq p_j\leq P_0\}$ with
$M_{0}$ and $P_{0}$ denoting upper bounds for $m$ and $p_{j}$,
respectively. The large sample behvavior of the minimum description
length criterion is given in the next theorem. Its proof can be found
in the \hyperref[sec:proof]{Appendix}.
%
%
\begin{theorem}\label{th:1}
Assume that the conditions of Proposition~\ref{prop:3} are satisfied
and let the number of break points $m^0$ be known. Then estimating the
piecewise quantile autoregressive model specified in (\ref{eqn:qar5})
at any single quantile $\tau\in(0,1)$ leads to
\[
\hat{\lambda}_{j}\rightarrow\lambda^{0}_{j}\qquad
\mbox{with probability one } (n\rightarrow\infty)
\]
for all $j=1,2,\ldots,m^0$, where $\hat{\lambda}=(\hat{\lambda
}_{1},\ldots,\hat{\lambda}_{m^0})$ is the minimizer of the criterion
function (\ref{eqn:qar10}).
\end{theorem}

The following \newtext{corollary} extends the result of Theorem~\ref
{th:1} to the multiple quantile case. Its verification is also provided
in the \newtext{\hyperref[sec:proof]{Appendix}}.
%
%
\begin{corollary}\label{cor:1}
Assume that the conditions of Proposition~\ref{prop:3} are satisfied.
Let the number of break points $m^0$ be known and assume that the break
locations as well as the autoregressive orders are the same across the
quantiles under consideration. Then estimating the piecewise quantile
autoregressive model specified in (\ref{eqn:qar5}) at the collection
of quantiles $(\tau_{1},\ldots,\tau_{L})\in(0,1)^L$ leads to
\[
\hat{\lambda}_{j}\rightarrow\lambda^{0}_{j}\qquad
\mbox{with probability one } (n\rightarrow\infty)
\]
for all $j=1,2,\ldots,m^0$, where $\hat{\lambda}=(\hat{\lambda
}_{1},\ldots,\hat{\lambda}_{m^0})$ is the minimizer of the criterion
function (\ref{eqn:qar11}).
\end{corollary}

We remark that in practice the assumption of known $m^0$ is often
unrealistic. However, it is substantially more difficult to establish
consistency in the general case of unknown $m^0$. Even in the simpler
univariate change-point frameworks, where independent variables are
grouped into segments of identical distributions, only special cases
such as normal distributions and exponential families have been
thoroughly investigated; for example, \cite{Lee.cb.97,Yao88} as well
as \cite{Aue:Lee:2011} for image segmentation. The
reason for this is that sharp tail estimates for maxima of certain
squared Gaussian processes are needed which do not hold for
distributions with thicker tails.

\section{Practical minimization using genetic algorithms}
\label{sec:ga}

Practical minimization of the minimum description length criteria (\ref
{eqn:qar10}) and (\ref{eqn:qar11}) is not a trivial task. We propose
using genetic algorithms to solve this minimization problem.

Genetic algorithms are a class of stochastic optimization techniques.
They are based on the idea of Darwin's theory of natural selection.
Typically a genetic algorithm begins with a random population of
possible solutions to the optimization problems. These solutions are
known as \emph{chromosomes} and often represented in vector form.
These chromosomes are allowed to evolve over time through the so-called
\emph{crossover} and \emph{mutation} operations. The hope is that the
evolution process would ultimately lead to a chromosome which
represents a good answer to the optimization problem. Successful
applications of genetic algorithms for solving various optimization
problems can be found, for examples, in \cite{Davis:1991}.

For a similar piecewise AR modeling minimization problem, Davis \textit{et al.} \cite{Davis:Lee:Rodriguez-Yam:2006}
developed a genetic algorithm for
approximating the minimizer. We modified their genetic algorithm to
solve the present minimization problem. For conciseness, we only
describe the major differences between the genetic algorithm for the
present piecewise quantile autoregressive model fitting problem and the
one from \cite{Davis:Lee:Rodriguez-Yam:2006}. We refer the reader to
\cite{Davis:Lee:Rodriguez-Yam:2006} for complete details.
\begin{ChromosomeRepresentation*} For the current problem of
detecting break points for a non-stationary time series at a specific
quantile $\tau$, a chromosome should contain information of all the
break points $k_{j}$ as well as the quantile autoregression orders
$p_{j}$ for any $\mathcal{F}\in\mathcal{M}$, where $\mathcal{M}$
denotes the whole class of piecewise quantile autoregressive models. We
express a chromosome as a vector of $n$ integers: a chromosome $\boldsymbol{\delta}=(\delta_{1},\ldots,\delta_{n})$ is of length $n$ with gene
values $\delta_{t}$ defined as
\[
\delta_{t}=\lleft\{ %
\begin{array} {l@{\qquad}l}
-1,& \mbox{if no
break point at time }t,
\\
p_{j},& \mbox{if $t=k_{j-1}$ and for the $j$th piece we
choose the $\QAR(p_j)$ model at quantile $\tau$.}
\end{array}
\rright.
\]
In practice, we impose an upper bound $P_{0}$ on the order $p_{j}$ of
each quantile autoregressive process. For our numerical work, we set
$P_{0}=20$. While the algorithm is running, we also impose the
following constraint on each $\delta$: in order to have enough
observations for parameter estimation, each piecewise quantile
autoregressive process is required to have a minimum length $m_p$,
which is chosen as a function of the order $p_j$ of the piecewise
process; their values are listed in Table~\ref{table1}.
\end{ChromosomeRepresentation*}
\begin{IslandModelandConvergence*} The Island Model was also
applied to speed up the convergence rate. We used 40 islands with
subpopulation size 40, performed a migration for every 5 generations,
and migrated 2 chromosomes during each migration.
And at the end of each migration the overall best chromosome that has
the smallest minimum description length value is selected. If this best
chromosome does not change for 20 consecutive migrations, or the total
number of generations exceeds 100, the genetic algorithm stops and the
best chromosome is taken as the solution to the optimization problem.
\end{IslandModelandConvergence*}
\ignore{
\textsc{Elitist step and Island Model:} In order to converge to an
optimal solution at a faster rate, two additional steps, the \emph
{elitist} step and the \emph{Island model}, are performed. In the
elitist step, the worst chromosome of the next generation is replaced
by the best chromosome of the current generation. This conserves the
best chromosome in each generation and guarantees the monotonicity in
the search process. In conducting the island model, parallel
implementations are applied and $NI$ (number of islands) genetic
algorithms are simultaneously performed in $NI$ different
sub-populations. After every $J$ generations, the worst $H$ chromosomes
from the $j$-th island are replaced by the best $H$ chromosomes from
the $(j-1)$-th island, $j=2,\ldots,NI$ (for $j=1$ the best $H$
chromosomes are migrated from the $NI$-th island). In our simulations
we set $NI=40$, $J=5$, $H=2$ and a sub-population size of 40.

\textsc{Convergence:} At the end of each migration the overall best
chromosome that has the smallest minimum description length value is
selected. If this best chromosome does not change for 20 consecutive
migrations, or the total number of generations exceeds 100, the genetic
algorithm stops and the best chromosome is taken as the solution to the
optimization problem.
}

%
%
\begin{table}[t]
\tablewidth=\textwidth
\tabcolsep=0pt
\caption{Values of $m_{p}$ used in the genetic algorithm}\label{table1}
\begin{tabular*}{\textwidth}{@{\extracolsep{\fill}}lllllllll@{}}
\hline
&\multicolumn{8}{l}{$p$}\\[-5pt]
&\multicolumn{8}{l}{\hrulefill}\\
& 0--1 & 2 & 3 & 4 & 5 &6 &7--10 &11--20 \\
\hline
$m_{p}$ & 10 & 12 & 14 & 16 & 18 & 20 &25 &50 \\ \hline
\end{tabular*}
\end{table}

\section{Simulation studies}
\label{sec:sim}

\subsection{Preliminaries}

In this section, four sets of simulation experiments are conducted to
evaluate the empirical performance of the proposed method for fitting
piecewise stationary quantile autoregressions. We shall compare the
results from our method with the Auto-PARM method of \cite
{Davis:Lee:Rodriguez-Yam:2006}, who developed an automatic procedure
for fitting piecewise autoregressive processes.
In each set of experiments, the results are based on 500 repetitions.
For the proposed method, we estimated the structural changes at
individual quantiles $\tau=0.25$, $0.5$ and $0.75$, as well as jointly
at $(0.25, 0.5, 0.75)$
using equal weights for the three quantiles. For convenience, we will
report the \emph{relative}\vspace*{2pt} locations of break points defined as $\hat
{\lambda}_{j}=\hat{k}_{j}/n$ for $j=1,\ldots,\hat{m}$.

\subsection{Piecewise $\AR(2)$ processes}

This simulation experiment is designed to compare the performance of
the proposed method and Auto-PARM in a linear autoregressive process setting
favoring the latter. The data generating process is
%
%
\begin{equation}
\label{eqn:sim1} y_{t}=\lleft\{ %
\begin{array} {l@{\qquad}l} 0.5y_{t-1}+0.3y_{t-2}+\varepsilon_{t}& (1
\leq t\leq n/2),
\\
-0.5y_{t-1}-0.7y_{t-2}+\varepsilon_{t}&
(n/2<t\leq3n/4),
\\
1.3y_{t-1}-0.5y_{t-2}+\varepsilon_{t}&
(3n/4<t\leq n),
\end{array}
\rright.
\end{equation}
where $(\varepsilon_{t})$ are independent standard normal, and
$n=1024$ and $2048$.
%
%
\begin{table}
\tablewidth=\textwidth
\tabcolsep=0pt
\caption{Summary of the estimated number of break points $\hat m$ for
the proposed procedure for the process (\protect\ref{eqn:sim1}) with $n=1024$.
Mean (standard deviation (Std)) of the relative break point location is
reported where applicable. If mult is specified for the quantile, it
refers to the multiple case $\tau=(0.25,0.50,0.75)$. The rows labeled
Auto-PARM give the results for that method}
\label{table:sim1}
\begin{tabular*}{\textwidth}{@{\extracolsep{\fill}}lllllll@{}}
\hline
&\multicolumn{6}{l}{$\hat m$}\\[-5pt]
&\multicolumn{6}{l}{\hrulefill}\\
& \multicolumn{1}{l}{0}&\multicolumn{2}{l}{1} & \multicolumn{2}{l}{2}
& \multicolumn{1}{l}{3} \\[-5pt]
& \multicolumn{1}{l}{\hrulefill}&\multicolumn{2}{l}{\hrulefill} &
\multicolumn{2}{l}{\hrulefill}
& \multicolumn{1}{l}{\hrulefill} \\
\multicolumn{1}{l}{$\tau$} & \multicolumn{1}{l}{\%} &
\multicolumn{1}{l}{\%} & \multicolumn{1}{l}{Mean (Std)} &
\multicolumn{1}{l}{\%} & \multicolumn{1}{l}{Mean (Std)} &
\multicolumn{1}{l}{\%} \\
\hline
$0.25$ & $1.2$ & $23.2$ & $0.759\ (0.016)$ & $75.6$ & $0.501\ (0.024)$ &
0 \\
& & & & & $0.747\ (0.012)$ & \\[3pt]
$0.50$ & $0$ & $\hphantom{2}3.6$ & $0.757\ (0.012)$ &
$96.4$ & $0.504\ (0.021)$ & 0 \\
& & & & & $0.747\ (0.011)$ & \\[3pt]
$0.75$ & $0.6$ & $19.8$ & $0.756\ (0.014)$ & $79.6$ & $0.501\ (0.025)$ &
0 \\
& & & & & $0.747\ (0.013)$ & \\[3pt]
mult & $0$ & $14.2$ & $0.750\ (0.013)$ & $85.8$ & $0.503\ (0.023)$ & 0 \\
& & & & & $0.748\ (0.012)$ & \\[3pt]
Auto-PARM & $0$ &$\hphantom{1}0$ & & $99.6$ & $0.501\ (0.004)$ &
$0.4$ \\
& & & & & $0.751\ (0.002)$ & \\ \hline
\end{tabular*}
\vspace*{-3pt}
\end{table}

For each simulated process we applied both procedures to locate the
break points. We recorded the number of break points detected by each
method, together with their relative locations. These numbers are
summarized in Tables~\ref{table:sim1} and \ref{table:sim1-2048}. From
Table~\ref{table:sim1}, we observe that, for the case $n=1024$, the
performance of Auto-PARM is slightly better than for the proposed
method at the median and is better at the other two quantiles under
consideration. However, as $n$ increased to 2048, the performance of
the quantile autoregression procedure improved and is comparable with
Auto-PARM both in terms of finding the correct number of breaks and
their locations, as can be seen from Table~\ref{table:sim1-2048}.

We have repeated the same experiment but with innovations distributed
as the $t$-distribution with 5 degrees of freedom. In this case, our
method outperformed Auto-PARM for all quantiles tested. Due to space
limitation, tabulated results are omitted.
%
%
\begin{table}[t]
\tablewidth=8cm
\tabcolsep=0pt
\caption{Similar to Table \protect\ref{table:sim1} except for
$n=2048$}\label{table:sim1-2048}
\begin{tabular*}{8cm}{@{\extracolsep{\fill}}llll@{}}
\hline
&\multicolumn{3}{l}{$\hat m$} \\[-5pt]
&\multicolumn{3}{l}{\hrulefill} \\
& \multicolumn{2}{l}{2} & \multicolumn{1}{l}{3} \\[-5pt]
& \multicolumn{2}{l}{\hrulefill} & \multicolumn{1}{l}{\hrulefill} \\
\multicolumn{1}{l}{$\tau$} & \multicolumn{1}{l}{\%} &\multicolumn
{1}{l}{ Mean (Std)} & \multicolumn{1}{l}{\%} \\
\hline
$0.25$ & $\hphantom{1}99.2$ & $0.503\ (0.015)$ &
$0.8$ \\
& & $0.747\ (0.008)$ & \\[3pt]
$0.50$ & $\hphantom{1}99.4$ & $0.503\ (0.012)$ &
$0.6$ \\
& & $0.744\ (0.006)$ & \\[3pt]
$0.75$ & $\hphantom{1}99.6$ & $0.503\ (0.015)$ &
$0.4$ \\
& & $0.748\ (0.007)$ & \\[3pt]
mult & $\hphantom{1}99.4$ & $0.504\ (0.013)$ & $0.6$ \\
& & $0.748\ (0.007)$ & \\[3pt]
Auto-PARM & $100$ & $0.501\ (0.002)$ & 0 \\
& & $0.750\ (0.001)$ & \\
\hline
\end{tabular*}
\end{table}

%
%
\begin{figure}[b]

\includegraphics{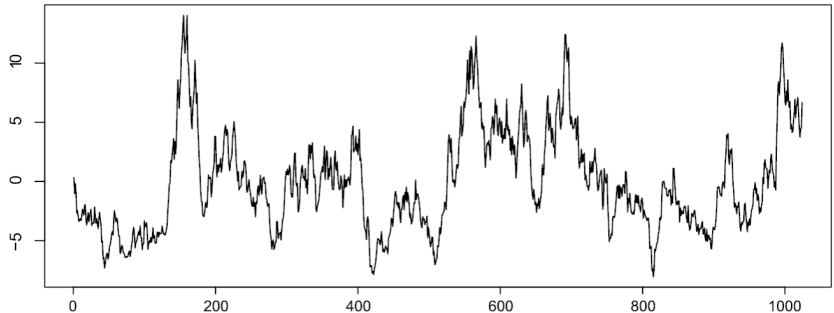}

\caption{A typical realization for the process in (\protect\ref{eqn:sim2}).}
\label{fig:sim2}
\end{figure}

\subsection{$\QAR(1)$ processes exhibiting explosive behavior}
\label{sec:sim2}

The data generating mechanism in this simulation follows the $\QAR(1)$ process
%
%
\begin{equation}
\label{eqn:sim2} y_{t}=(0.85+0.25u_t)y_{t-1}+
\Phi^{-1}(u_t),
\end{equation}
where $(u_t)$ is a sequence of independent standard uniform random
variables and $\Phi$ the standard normal distribution function. Shown
in Figure {\ref{fig:sim2}} is a typical realization. There is no
structural break in this series but from the plot one can see that it
exhibits explosive behavior in the upper tail. Processes such as this
one seem to be capable of modeling certain macroeconomic time series;
for example, interest rate data. We will
revisit this issue in Section~\ref{sec:app} below. While our method
does not detect break points at any of the quantiles tested, only about
one-third of the results from Auto-PARM lead to the correct conclusion;
the numbers of break points detected by their method are summarized in
Table~\ref{table:sim2}. It is apparent that it is much less tolerant
to nonlinearity.

\subsection{Piecewise $\AR(1)$ processes with changes in certain
quantile ranges}

In this simulation experiment, the nonstationary time series is
generated from the model
%
%
\begin{equation}
\label{eqn:sim3} y_{t}=\lleft\{ %
\begin{array} {l@{\qquad}l}
\bigl\{0.5I(\tau\leq0.2)+0.8I(\tau>0.2)\bigr\} y_{t-1}+
\varepsilon_{t}& (1\leq t\leq n/2),
\\
0.5y_{t-1}+\varepsilon_{t}& (n/2<t\leq n),
\end{array}
\rright.
\end{equation}
where $(\varepsilon_t)$ are independent asymmetric Laplace with
parameter $0.4$ for $t\leq n/2$ and independent asymmetric
Laplace with parameter $0.6$ for $t>n/2$.

%
%
\begin{table}
\tablewidth=\textwidth
\tabcolsep=0pt
\caption{Relative frequencies of the number of break points estimated
from Auto-PARM for the process (\protect\ref{eqn:sim2}) with $n=1024$.
Independent of the specific quantile it was applied to, the proposed
methodology always correctly chose $\hat m=0$}\label{table:sim2}
\begin{tabular*}{\textwidth}{@{\extracolsep{4in minus 4in}}lllllll@{}}
\hline
&\multicolumn{6}{l}{Number of break points} \\[-5pt]
&\multicolumn{6}{l}{\hrulefill} \\
& 0 & 1 & 2 & 3 & 4 &5 \\
\hline
Relative frequency &33.8 &35.2 &23.8 &
5.6 & 1.4 & 0.2 \\ \hline
\end{tabular*}
\end{table}

%
%
\begin{table}[b]
\tablewidth=11cm
\tabcolsep=0pt
\caption{Similar to Table \protect\ref{table:sim1} except for the
process (\protect\ref{eqn:sim3}) with $n=1024$}\label{table:sim3}
\begin{tabular*}{11cm}{@{\extracolsep{\fill}}llllll@{}}
\hline
&\multicolumn{5}{l}{$\hat m$} \\[-5pt]
&\multicolumn{5}{l}{\hrulefill} \\
& \multicolumn{1}{l}{0} & \multicolumn{2}{l}{1} & \multicolumn
{1}{l}{2} & \multicolumn{1}{l}{3} \\[-5pt]
& \multicolumn{1}{l}{\hrulefill} & \multicolumn{2}{l}{\hrulefill} &
\multicolumn{1}{l}{\hrulefill} & \multicolumn{1}{l}{\hrulefill} \\
\multicolumn{1}{l}{$\tau$} & \multicolumn{1}{l}{\%} & \multicolumn
{1}{l}{\%} & \multicolumn{1}{l}{Mean (Std)} &
\multicolumn{1}{l}{\%} & \multicolumn{1}{l}{\%} \\
\hline
$0.25$ & $83.4$ & $16.6$ & $0.527\ (0.096)$ & 0 & 0 \\
$0.50$ & $\hphantom{8}1.5$ & $98.5$ & $0.503\ (0.038)$ & 0 & 0 \\
$0.75$ & $24.4$ & $75.6$ & $0.479\ (0.055)$ & 0 & 0 \\
mult & $35.2$ & $64.8$ & $0.498\ (0.046)$ &
0 & 0 \\
Auto-PARM & $51.0$ & $44.4$ & $0.487\ (0.181)$ & $4.0$ & $0.6$ \\ \hline
\end{tabular*}
\end{table}

For this process, results from our method and Auto-PARM are reported in
Table~\ref{table:sim3} in a similar manner as in Table~\ref
{table:sim1}. Not reported in this table is the fact that, when the
coefficients of $y_{t-1}$ in the two pieces are the same (which happens
for quantiles $\tau\leq0.2$), then the proposed procedure does
not detect any break points even though the residuals of the two pieces
are slightly different. For the quantile at $\tau=0.25$ which
is close to the threshold at which the autoregressive coefficient
changes, our method detected a (nonexisting) break point in 16\% of the
simulation runs. On the other hand, when $\tau\geq0.5$, the
quantile autoregression method performs reasonably well, especially at
the median where the performance is excellent. Also at $\tau=0.75$ it
outperforms Auto-PARM. When estimating jointly at $\tau
=(0.25, 0.5, 0.75)$, the percentage of detecting
the correct number of break points is not as high as at $\tau=0.5$ due
to the inclusion of the quantiles at $\tau=0.25$ and
$\tau=0.75$, indicating that care has to be exercised if
quantiles are jointly specified. We can also see that the performance
of our method is better than that of Auto-PARM in both percentage and
accuracy (in terms of smaller standard deviations) for this simulation
example. In Table~\ref{table:sim3order}, we summarize the proposed
procedure's estimates of the quantile autoregression orders for the
above process at $\tau=0.5$, and we can see that most of the
segments are correctly modeled as $\QAR(1)$ processes.
%
%
\begin{table}[t]
\tablewidth=200pt
\tabcolsep=0pt
\caption{Relative frequencies of the quantile autoregression orders
selected by the proposed method at $\tau=0.5$ for the
realizations from the process (\protect\ref{eqn:sim3})}\label
{table:sim3order}
\begin{tabular*}{200pt}{@{\extracolsep{\fill}}llllll@{}}
\hline
&\multicolumn{5}{l}{Order}\\[-5pt]
&\multicolumn{5}{l}{\hrulefill}\\
& \multicolumn{1}{l}{1} & \multicolumn{1}{l}{2} & \multicolumn{1}{l}{3}
& \multicolumn{1}{l}{4} & \multicolumn{1}{l}{5} \\ \hline
$p_{1}$ & $80.3$ & $15.7$ & $2.6$ & $1.4$ &
0 \\
$p_{2}$ & $72.4$ & $19.2$ & $6.6$ & $1.4$ &
$0.4$ \\
\hline
\end{tabular*}
\end{table}

\subsection{Higher-order QAR processes}

\begin{figure}

\includegraphics{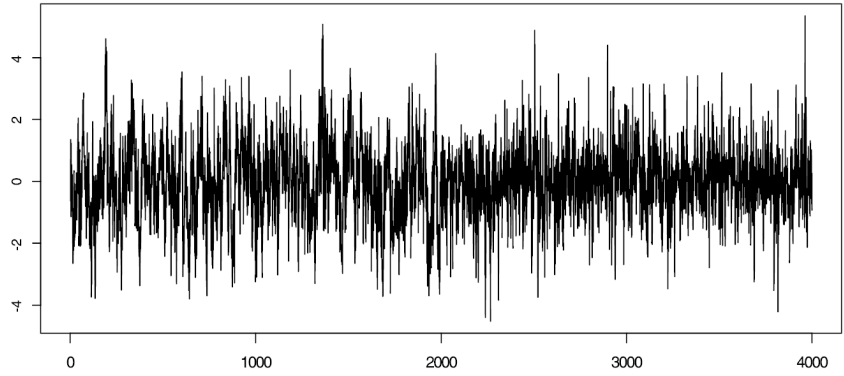}

\caption{A typical realization for the process in (\protect\ref{eqn:sim4}).}\vspace*{10pt}
\label{fig:sim4}
\end{figure}

In this experiment, the data generating process is
%
%
\begin{equation}
\label{eqn:sim4} y_{t}=\lleft\{ %
\begin{array} {l@{\qquad}l}
(0.2+0.1u_t)y_{t-1}+(0.5+0.1u_t)y_{t-2}+\epsilon_{t} & (1\leq t\leq n/2),
\\
0.7u_ty_{t-1}+\epsilon_{t}& (n/2<t\leq
n),
\end{array}
\rright.
\end{equation}
where $(u_t)$ is a sequence of independent standard uniform random
variables, $(\epsilon_t)$ are independent standard normal for $t \leq
n/2$, and independent asymmetric Laplace with parameter 1 for $t >
n/2$. A typical realization is displayed in Figure~\ref{fig:sim4}, and
break detection results from our method for this process are reported
in Table~\ref{table:sim4}. One can see that our method has
successfully detected one break with very high probability in most
considered cases, and that the detected relative locations are also
very close to the true location.

In order to assess the performance of the MDL criterion for order
selection in $\QAR(p)$ models for $p>1$, we tabulated the relative
frequencies of the order selected by the proposed method for the first
piece of process (\ref{eqn:sim4}) in Table~\ref{table:sim4_2}. The
proposed method never underestimates the order, but only achieves about
50\% accuracy. At first sight, these correct estimation rates seem to
be relatively low. However, in the break point detection context, the
problem of order estimation seems to be hard even for linear AR
processes (of higher order), as is seen in Table~3 of \cite
{Davis:Lee:Rodriguez-Yam:2006}, where Auto-PARM only gave around 65\%
correct estimation rates for $\AR(2)$ processes. Thus, we believe that a
50\% correct rate is not unreasonable for $\QAR(p)$ models.
%
%
\begin{table}
\tablewidth=250pt
\tabcolsep=0pt
\caption{Similar to Table \protect\ref{table:sim1} except for the
process (\protect\ref{eqn:sim4}) with $n=4000$}
\label{table:sim4}
\begin{tabular*}{250pt}{@{\extracolsep{\fill}}lllll@{}}
\hline
& \multicolumn{4}{l}{$\hat m$} \\[-5pt]
& \multicolumn{4}{l}{\hrulefill} \\
&\multicolumn{1}{l}{0} &\multicolumn{2}{l}{1} & \multicolumn
{1}{l}{2} \\[-5pt]
&\multicolumn{1}{l}{\hrulefill} &\multicolumn{2}{l}{\hrulefill} &
\multicolumn{1}{l}{\hrulefill} \\
\multicolumn{1}{l}{$\tau$} &\multicolumn{1}{l}{\%}& \multicolumn
{1}{l}{\%} & \multicolumn{1}{l}{Mean (Std)} & \multicolumn{1}{l}{\%}
\\
\hline
$0.25$ &4.0& \hphantom{1}95.5 & 0.517\ (0.049)&0.5\\
$0.50$ &0& \hphantom{1}98.5 & 0.505\ (0.039)&1.5 \\
$0.75$ &3.0& \hphantom{1}97.0 & 0.508\ (0.052) &0 \\
mult &0& 100.0 & 0.509\ (0.045) & 0.5\\
\hline
\end{tabular*}
%
\vspace*{30pt}
%
%
\tablewidth=\textwidth
\tabcolsep=0pt
\caption{Relative frequencies of the quantile autoregression orders
selected by the proposed method at different $\tau$ values ($\tau=$
0.25, 0.50, 0.75, and mult) for the first piece in the process
(\protect\ref
{eqn:sim4}). The true order is 2}
\label{table:sim4_2}
\begin{tabular*}{\textwidth}{@{\extracolsep{\fill}}llllllll@{}}
\hline
$\tau$ & 1&2 & 3 & 4 & 5 & 6 & $\geq7$ \\ \hline
0.25 & 0& 48.69& 31.41&15.71& 2.09& 1.57& 0.52 \\
0.50 & 0& 51.78& 26.40& 12.18& 5.58 & 2.03& 2.03 \\
0.75 & 0&55.15&22.68& 11.86& 7.73& 1.55& 1.05 \\
mult & 0 &50.50& 26.00& 14.50& 5.00& 2.00& 2.00 \\ \hline
\end{tabular*}
\end{table}

\ignore{
%
%
\begin{table}[!ht]
\begin{center}
\begin{tabular}{|c|c|@{ }ccccccccc|}\hline
$\tau$ & Order& 1&2 & 3 & 4 & 5 &6&7&8&9 \\ \hline
0.25 & $p_{1}$ & 0& 48.69& 31.41&15.71& 2.09& 1.57& 0& 0.52 &0 \\
& $p_{2}$ & 46.07&31.94& 12.57& 5.24& 2.62& 1.05 & 0& 0.52& 0 \\ \hline
0.50 & $p_1$ & 0& 51.78& 26.40& 12.18& 5.58 & 2.03& 2.03& 0&0 \\
& $p_2$ & 50.25& 24.37& 13.20& 5.58& 4.06& 1.02& 1.02&1.02&0 \\ \hline
0.75 & $p_1$ & 0&55.15&22.68& 11.86& 7.73& 1.55& 0.52&0&0.52 \\
& $p_2$ & 48.97& 29.90& 10.31&7.22& 2.06& 1.03& 0.52&0&0 \\ \hline
mult & $p_1$ & 0 &50.50& 26.00& 14.50& 5.00& 2.00& 2.00& 0& 0 \\
& $p_2$ &40.00& 31.50& 16.00& 7.50& 2.50& 2.00& 0.5&0&0 \\ \hline
\end{tabular}

\caption{Relative frequencies of the quantile autoregression orders
selected by the proposed method at different $\tau$ values ($\tau$ =
0.25, 0.50, 0.75, and mult).}\label{table:sim4_2}
\end{center}
\end{table}
}

\subsection{Stochastic volatility models}
\label{sec:svm}

The simulation section concludes with an application of the proposed
methodology to stochastic volatility models (SVM) often used to fit
financial time series; see \cite{Shephard:Andersen:2009} for a recent
overview. It should be noted that the proposed quantile methodology and
Auto-PARM are not designed to deal with this type of model as it
consists of uncorrelated random variables exhibiting dependence in
higher-order moments. However, SVM are used to compare the two on a
data generating process different from nonlinear QAR and linear AR time
series. Following Section~4.2 of \cite{Davis:Lee:Rodriguez-Yam:2008},
the process
%
%
\begin{equation}
\label{eq:svm1} y_t = \sigma_t\xi_t=e^{\alpha_t/2}
\xi_t,
\end{equation}
is considered, where $\alpha_t=\gamma+\phi\alpha_{t-1}+\eta_t$.
The following two-piece segmentations were compared:
\begin{eqnarray*}
&&\mbox{Scenario A} \quad\mbox{Piece 1:}\quad \gamma
=-0.8106703,\qquad \phi=0.90,\qquad (
\eta_t)\sim\mbox{ i.i.d. } N(0,0.45560010),
\\
&&\hphantom{\mbox{Scenario A} \quad}\mbox{Piece 2:}\quad \gamma
=-0.3738736,\qquad \phi=0.95, \qquad(\eta_t)\sim\mbox{
i.i.d. } N(0,0.06758185),
\end{eqnarray*}
while $(\xi_t)\sim\mbox{ i.i.d. } N(0,1)$ for both pieces, and
\begin{eqnarray*}
&&\mbox{Scenario B}\quad \mbox{Piece 1:}\quad \gamma=-0.8106703,
\qquad\phi=0, \qquad(
\xi_t)\sim\mbox{ i.i.d. } N(0,1),
\\
&&\hphantom{\mbox{Scenario B}\quad}\mbox{Piece 2:}\quad \gamma
=-0.3738736,\qquad \phi=0,\qquad (\xi_t)\sim\mbox{
i.i.d. } N(0,4),
\end{eqnarray*}
while $(\eta_t)\sim\mbox{ i.i.d. } N(0,0.5)$ for both pieces.
Scenario A corresponds to a change in dynamics of the volatility
function $\sigma_t$, Scenario B basically to a scale change.

\begin{figure}

\includegraphics{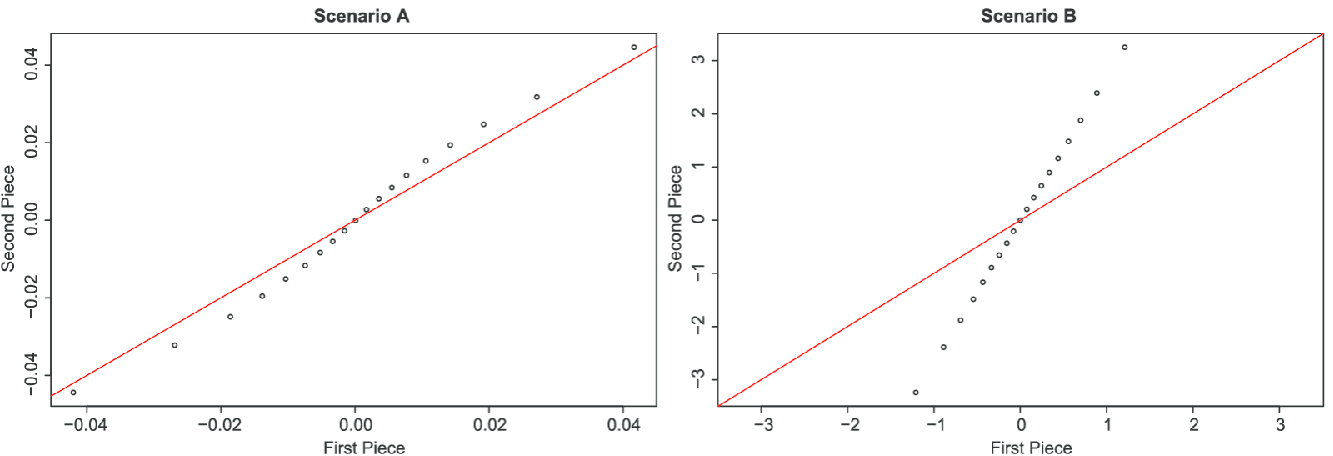}

\caption{Empirical quantile--quantile plot for the SVM process
specified under Scenario A (left panel) and Scenario B (right panel).
The $x$-axis ($y$-axis) shows the empirical quantiles of Piece 1 (Piece
2). The 45 degree line is given for ease of comparison.}
\label{fig:quantile-quantile}
\end{figure}

Scenario A was considered in \cite{Davis:Lee:Rodriguez-Yam:2008}.
These authors developed a method tailored to deal with financial time
series of SVM and GARCH type. The method, termed Auto-Seg, was able to
detect one break in 81.8\% of 500 simulation runs and detected no break
otherwise. On this data, Auto-PARM tends to use a too fine segmentation
as 62.4\% of the simulations runs resulted in two or more estimated
break points. One (no) breakpoint was detected in 21.2\% (16.4\%) of
the cases. The proposed method failed to detect any changes at any of
the tested quantiles ($\tau=0.05, 0.10, 0.25, 0.50, 0.75, 0.90,
0.95$). It should be noted, however, that there is no change at the
median and changes in the other quantiles are very hard to find as is
evidenced by Figure~\ref{fig:quantile-quantile}, which displays the
averaged (over 50 simulation runs) empirical quantile--quantile plot
from the first and the second segment of the two-piece Scenario A process.
%
%

The results for Scenario B are summarized in Table~\ref{table:svm}. It
can be seen that, for the proposed method, the scale change, is
detected at the more extreme quantiles ($\tau=0.05, 0.10, 0.90, 0.95$)
with very good accuracy and with reasonable accuracy at intermediate
quantiles ($\tau=0.25$ and $\tau=0.75$), while no change is found
(correctly) at the median $\tau=0.50$, reflecting that the proposed
procedure describes the local behavior of the SVM process adequately.
Auto-PARM does the same on a global basis.
%
%
\begin{table}[b]
\tablewidth=\textwidth
\tabcolsep=0pt
\caption{Summary of the estimated number of break points $\hat m$ for
the proposed procedure and Auto-PARM for the process \protect\eqref{eq:svm1}
with specifications given under Scenario B}
\label{table:svm}
\begin{tabular*}{\textwidth}{@{\extracolsep{\fill}}lllllllll@{}}
\hline
& \multicolumn{7}{l}{$\tau$}& \\[-5pt]
& \multicolumn{7}{l}{\hrulefill}& \\
\multicolumn{1}{l}{$\hat m$}&  \multicolumn{1}{l}{0.05} &
\multicolumn{1}{l}{0.10} & \multicolumn{1}{l}{0.25}
& \multicolumn{1}{l}{0.50} & \multicolumn{1}{l}{0.75} & \multicolumn
{1}{l}{0.90} & \multicolumn{1}{l}{0.95} &
\multicolumn{1}{l}{Auto-PARM} \\ \hline
0 &   \hphantom{9}0.4\% & \hphantom{9}0.2\% & 32.6\% & 100.0\% & 29.6\% &
\hphantom{10}0.0\% & \hphantom{9}0.6\% & \hphantom{9}0.2\% \\
1 & 99.6\% & 99.8\% & 67.4\% & \hphantom{10}0.0\% & 70.4\% & 100.0\% & 99.4\% &
99.6\% \\
2 &  \hphantom{9}0.0\% & \hphantom{9}0.0\% & \hphantom{9}0.0\% &
\hphantom{10}0.0\% &
\hphantom{9}0.0\% & \hphantom{10}0.0\% & \hphantom{9}0.0\% & \hphantom
{9}0.2\% \\ \hline
\end{tabular*}
\end{table}

\section{Real data applications}
\label{sec:app}

\subsection{Treasury bill data}

\begin{figure}

\includegraphics{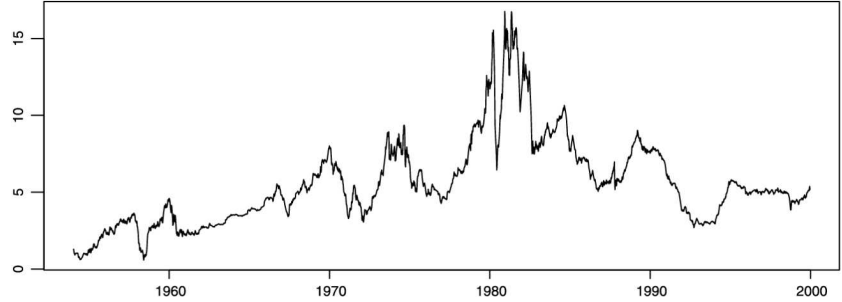}

\caption{Three-month treasury bills (01/1954 to 12/1999).} \label{fig:bills03}
\end{figure}

Treasury bills are short-term risk-free investments that are frequently
utilized by investors to hedge portfolio risks. In this application,
the observations are three-month treasury bills from the secondary
market rates in the United States, ranging from January 1954 to
December 1999. The weekly data can be found at the website \url
{http://research.stlouisfed.org/fred2/series/TB3MS} and are displayed
in Figure~\ref{fig:bills03}.
%
%

It can be seen from Figure~\ref{fig:bills03} that the time series
exhibits obvious explosive behavior in the upper tail. In many
instances similar time series would be viewed as random walks and
sophisticated testing procedures would have to be applied to either
confirm or reject what is known as unit-root hypothesis; see, for
example, \cite{Paparoditis:Politis:2003,Paparoditis:Politis:2005} for
more. As in Section~\ref{sec:sim2}, Auto-PARM aims in this case at
partitioning the series into segments with structures mimicking linear
behavior. In the present case, this leads to 15 segments. On the other
hand the proposed procedure does not detect break points at any of the
quantiles tested ($\tau=0.05, 0.10,\ldots,0.90,0.95$), thus
indicating that with the use of some extra
parameters a more parsimonious stationary but nonlinear modeling is
possible for this data set. Using a $\QAR(2)$ model with cubic polynomial
coefficients in the uniform random variables $(u_t)$, the data can be
approximated via the following model with 12 parameters:
%
%
\begin{equation}
\label{eqn:bills03} y_t=\theta_0(u_t)+
\theta_1(u_t)y_{t-1}+\theta_2(u_t)y_{t-2},
\end{equation}
where
\begin{eqnarray*}
\theta_0(u_t)&=&-0.0144+0.2264u_t-0.5448u_t^2+0.3848u_t^3,
\\
\theta_1(u_t)&=&1.3721-0.9635u_t+1.5312u_t^2-0.6939u_t^3,
\\
\theta_2(u_t)&=&-0.4394+1.3154u_t-2.1945u_t^2+1.1353u_t^3.
\end{eqnarray*}

\begin{figure}

\includegraphics{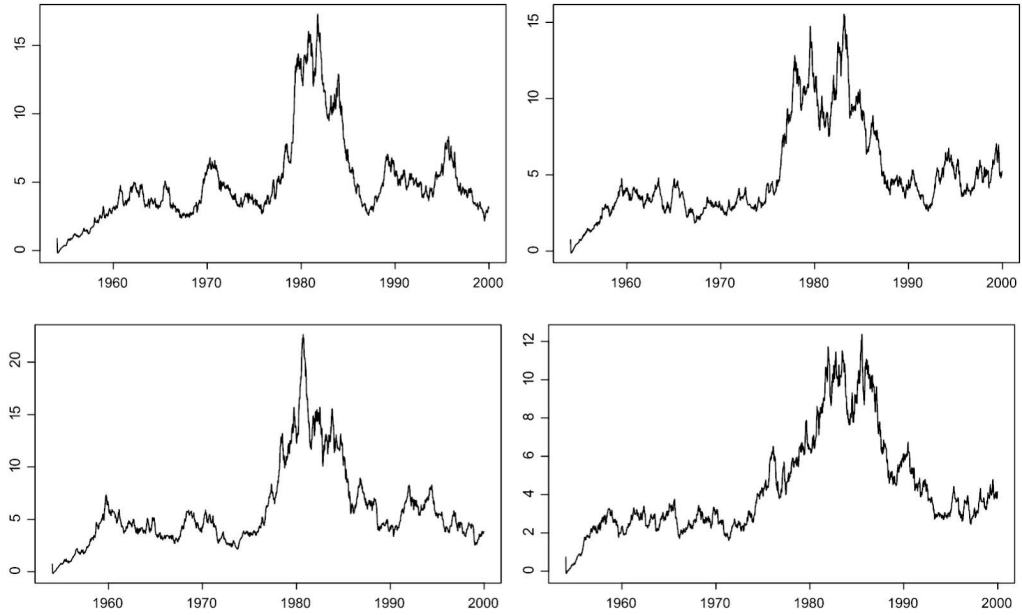}

\caption{Four typical realizations of the process in (\protect\ref
{eqn:bills03}).} \label{fig:generated}
\end{figure}

Figure~\ref{fig:generated} depicts several realizations generated by
the estimated model (\ref{eqn:bills03}), which all show a pattern
closely resembling the data in Figure~\ref{fig:bills03}. This example
illustrates that quantile autoregressions can expand the modeling
options available to the applied statistician as it accurately captures
temporary explosive behavior and nonlinearity.
%
%

\subsection{Monthly minimum temperature data}

\begin{figure}[b]

\includegraphics{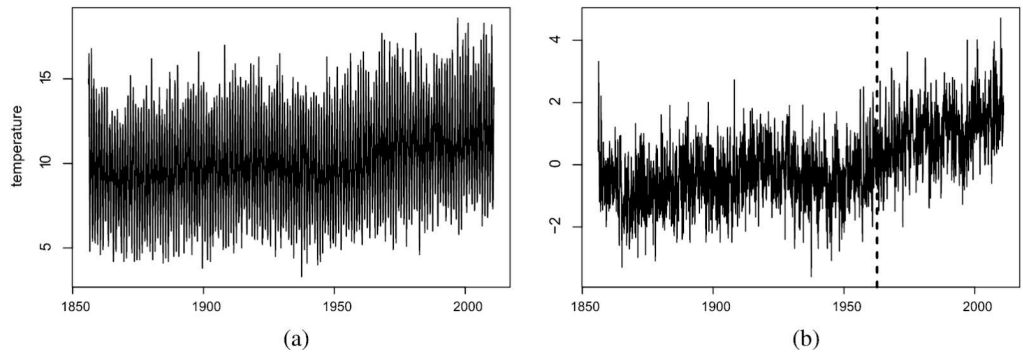}

\caption{(a) Monthly minimum air temperature in Melbourne, Australia
from January 1856 to December 2010. (b) Deseasonalized series. The
dashed line represents the estimated break point in August 1962.}
\label{fig:australia}
\end{figure}

In this section the monthly mean minimum temperature at Melbourne in
Australia is considered. The data set is obtainable from the Bureau
of Meteorology of the Australian Government (\href
{http://www.bom.gov.au/climate/data/}{http://  www.bom.gov.au/climate/data/}). The plots for the original
series and its deseasonalized version are shown in Figure~\ref
{fig:australia}. This data set has been investigated by \cite
{Koenker:2005} who pointed out that, due to the quantile dependent
behavior visible in the scatter plots, linear autoregressive models
are insufficient to describe the data. Our method was
applied to this data set at various quantiles and for all cases one
break point was found near the year 1960. This agrees with a visual
inspection of Figure~\ref{fig:australia}.
%
%
%
%
\begin{table}
\tablewidth=\textwidth
\tabcolsep=0pt
\caption{Estimated break points at different quantiles for the
Australian temperature data}\label{table:australia}
\begin{tabular*}{\textwidth}{@{\extracolsep{\fill}}lllll@{}}
\hline
&\multicolumn{4}{l}{Quantiles}\\[-5pt]
&\multicolumn{4}{l}{\hrulefill}\\
& \multicolumn{1}{l}{0.25} & \multicolumn{1}{l}{0.5} & \multicolumn
{1}{l}{0.75} & \multicolumn{1}{l}{mult} \\
\hline
Estimated break point & December 1960 &August 1963 & December 1958 &
August 1962 \\
\hline
\end{tabular*}
\end{table}

It can be seen from Table~\ref{table:australia} that the break point
location estimated with the multiple quantile procedure, set up with
equal weights for the three quantiles under consideration, is between
the break point locations estimated at the individual quantiles. This
should always be the case, as the requirement of simultaneous
occurrence of breaks automatically leads to a weighted average
interpretation. In general, one would ideally find weights that prefer
quantiles which stronger exhibit the structural break and attenuate the
impact of quantiles that are only marginally subjected to the break.
This would mean to more closely evaluate properties of the (piecewise)
density and distribution function of the underlying random process.

\section{Conclusions}
\label{sec:sum}

This article proposes a new segmentation procedure that helps breaking
down a given nonstationary time series into a number of stationary
pieces by means of quantile autoregression modeling. In contrast to
most of the existing literature, this is done either for individual
quantiles or across a collection of quantiles. The proposed method
utilizes the minimum description length principle and a genetic
algorithm to obtain the best segmentation. It has been proved that this
method is asymptotically consistent, and simulation results have
demonstrated that the finite sample performance of the proposed
procedure is quite good. Data applications are also provided with
satisfactory results. It can be seen in particular that our method can
add to second-order time series modeling by enriching the
statistician's tool box via the inclusion of nonlinearity, asymmetry,
local persistence and other distributional aspects. An interesting
problem for future research that shows some potential is the
investigation of the properties of the multiple quantile segmentation
procedure for the case of quantile-dependent break point locations,
thereby loosening the assumption of simultaneous breaks utilized in
this paper.


\begin{appendix}

\section*{Appendix: Proofs}
\label{sec:proof}
%
%
\begin{lemma}\label{prop:A.1}
If $(y_t\dvt t\in\mathbb{Z})$ follow a stationary $\QAR(p)$ model
such that the assumptions of Proposition~\ref{prop:3} are satisfied,
then with probability one and for all $\tau\in(0,1)$,
\[
\frac{1}n\sum_{t=1}^n
\rho_\tau(\hat\varepsilon_t)\to E\bigl\{
\rho_\tau(\varepsilon_1)\bigr\}\qquad (n\to\infty),
\]
where $\rho_\tau$ is the check function defined below (\ref{eqn:qar3}).
\end{lemma}
\begin{pf}
The assertion follows as in the proof of Lemma A.1 in \cite
{Aue:Cheung:Lee:Zhong:2014}.
\end{pf}
\ignore{Fix $\tau\in(0,1)$. The quantile autoregrssion equations
imply that $\hat\varepsilon_t=\varepsilon_t+X_t^\prime(\hat\theta
(\tau)-\theta(\tau))$. It follows from \cite{Koenker:2005} that
$\hat\theta(\tau)$ is strongly consistent for $\theta(\tau)$.
Therefore, with probability one,
\[
\frac{1}n |\sum_{t=1}^n(\hat
\varepsilon_t-\varepsilon_t) |=\frac{1}n |\sum
_{t=1}^nX_t^\prime
\bigl\{\hat\theta(\tau)-\theta(\tau)\bigr\} |\to0.
\]
Consequently, $\frac{1}n\sum_{t=1}^n\varepsilon_t\to E(\varepsilon
_1)$ with probability one by the strong law of large numbers. Since
$\rho_\tau(\cdot)$ is a continuous and measurable function, the
latter limit relation is also true if $\varepsilon_t$ is replaced by
$\rho_\tau(\epsilon_t)$.
}
%
%
\begin{lemma}\label{prop:A.2}
Let $(y_t\dvt t\in\mathbb{Z})$ be a piecewise stationary $\QAR(p)$
model that satisfies the assumptions of Proposition~\ref{prop:3} on
each of the segments. Let $\lambda^0=(\lambda^0_1,\ldots,\lambda
^0_{m^0})$ denote the true segmentation and choose $K=\lfloor\kappa
n\rfloor$, $M=\lfloor\mu n\rfloor$ with $0\leq\kappa<\mu\leq1$.
Then, with probability one for all $\tau\in(0,1)$,
\[
\frac{1}{M-K}\sum_{t=K+1}^{M}
\rho_\tau(\hat\varepsilon_t)\to L_\tau(\kappa,
\mu).
\]
The limit $L_\tau(\kappa,\mu)$ is the sum of two components, $A_\tau
(\kappa,\mu)$ and $B_\tau(\kappa,\mu)$, both of which are given in
the proof.
\end{lemma}
\begin{pf}
There are two cases to consider, namely (1) $K$ and $M$ are contained
in the same segment and (2) $K$ and $M$ are in different segments.

For the case (1), Lemma~\ref{prop:A.1} implies immediately that
\[
\frac{1}{M-K}\sum_{t=K+1}^M
\rho_\tau(\hat\varepsilon_t)\to\rho_{\tau,j}=
A_\tau(\kappa,\mu).
\]
With $B_\tau(\kappa,\mu)=0$, one can set $L_\tau(\kappa,\mu
)=A_\tau(\kappa,\mu)$ and the limit is determined.

For the case (2), there are $1\leq j<J\leq m^0+1$ such that $\kappa\in
[\lambda_{j-1}^0,\lambda_j^0)$ and $\mu\in(\lambda_{J-1}^0,\lambda
_J^0]$. In addition to the residuals $\hat\varepsilon_t$ obtained
from fitting a QAR model to the observations $y_{K+1},\ldots,y_M$, one
also defines residuals $\hat\varepsilon_{t,\ell}$ obtained from
fitting a QAR model on the $\ell$th underlying (true) segment. If now
$t\in\{k_{\ell-1}^0+1,\ldots,k_\ell^0\}$ with $k_\ell^0=\lfloor
\lambda_\ell^0 n\rfloor$, then one gets the decomposition $\rho
_\tau(\hat\varepsilon_t)=\{\rho_\tau(\hat\varepsilon_t)-\rho
_\tau(\hat\varepsilon_{t,\ell})\}+\rho_\tau(\hat\varepsilon
_{t,\ell})$. The sum over the first terms on the right-hand side leads
to a positive bias term $B_\tau(\kappa,\mu)$ determined by the
almost sure limit relation
\begin{eqnarray*}
&&\frac{1}{M-K} \Biggl[\sum_{t=K+1}^{k_{j}^0}
\bigl\{\rho_\tau(\hat\varepsilon_t)-\rho_\tau(
\hat\varepsilon_{t,j})\bigr\} \\
&&{\hphantom{\frac{1}{M-K} \Biggl[}}+\sum_{\ell=j+1}^{J-1}
\sum_{t=k_{\ell-1}+1}^{k_\ell^0}\bigl\{\rho_\tau(
\hat\varepsilon_t)-\rho_\tau(\hat\varepsilon_{t,\ell})
\bigr\} +\sum_{t=k_{J-1}^0+1}^M\bigl\{
\rho_\tau(\hat\varepsilon_t)-\rho_\tau(\hat
\varepsilon_{t,J})\bigr\} \Biggr]
\\
&&\quad\to B_\tau(\kappa,\mu).
\end{eqnarray*}
The remaining segment residuals $\hat\varepsilon_{t,\ell}$ allow for
an application of Lemma~\ref{prop:A.1} to each of the underlying
(true) segments, so that, with probability one,
\begin{eqnarray*}
&&\frac{1}{M-K} \Biggl\{\sum_{t=K+1}^{k^0_j}
\rho_\tau(\hat\varepsilon_{t,j})+\sum
_{\ell=j+1}^{J-1}\sum_{t=k^0_{\ell
-1}+1}^{k^0_\ell}
\rho_\tau(\hat\varepsilon_{t,\ell})+\sum
_{t=k_{J-1}^0+1}^M\rho_\tau(\hat
\varepsilon_{t,J}) \Biggr\}
\\
&&\quad\to\frac{1}{\mu-\kappa} \Biggl\{\bigl(\lambda_j^0-
\kappa\bigr)\rho_{\tau
,j}+\sum_{\ell=j+1}^{J-1}
\bigl(\lambda_\ell^0-\lambda_{\ell-1}^0
\bigr)\rho_{\tau,\ell}+\bigl(\mu-\lambda_{J-1}^0\bigr)
\rho_{\tau,J} \Biggr\}
\\
&&\quad=A_\tau(\kappa,\mu),
\end{eqnarray*}
where $\rho_{\tau,j}=E\{\rho_\tau(\varepsilon_{k_j^0})\}$. Setting
$L_\tau(\kappa,\mu)=A_\tau(\kappa,\mu)+B_\tau(\kappa,\mu)$
completes the proof.
\end{pf}
\begin{pf*}{Proof of Theorem~\ref{th:1}}
Denote by $\hat\lambda=(\hat\lambda_1,\ldots,\hat\lambda_{m^0})$
and $\lambda^0=(\lambda^0_1,\ldots,\lambda^0_{m^0})$ the
segmentation chosen by the minimum description length criterion (\ref
{eqn:qar10}) and the true segmentation, respectively. The proof is
obtained from a contradiction argument, assuming that $\hat\lambda$
does not converge almost surely to $\lambda^0$. If that was the case,
then the boundedness of $\hat\lambda$ would imply that, almost surely
along a subsequence, $\hat\lambda\to\lambda^*=(\lambda_1^*,\ldots
,\lambda_{m^0}^*)$ as $n\to\infty$, where $\lambda^*$ is different
from $\lambda^0$. Two cases for neighboring $\lambda_{j-1}^*$ and
$\lambda_j^*$ have to be considered, namely (1) $\lambda_{j^\prime
}^0\leq\lambda_{j-1}^*<\lambda_j^*\leq\lambda_{j^\prime}^0$ and
(2) $\lambda_{j^\prime-1}^0\leq\lambda_{j-1}^*<\lambda_{j^\prime
}^0<\cdots<\lambda_{j^\prime+J}^0<\lambda_j^*\leq\lambda
_{j^\prime+J+1}^0$ for some positive integer $J$.

For the case (1), Lemma~\ref{prop:A.1} implies that, almost surely,
\[
\lim_{n\to\infty}\frac{1}n\sum_{t=\hat k_{j-1}+1}^{\hat k_j}
\rho_\tau(\hat\varepsilon_t) \geq\bigl(
\lambda_j^*-\lambda_{j-1}^*\bigr)\rho_{\tau,j^\prime},
\]
where $\rho_{\tau,j^\prime}=E\{\rho_\tau(\varepsilon_{k_{j^\prime
}^0})\}$.
For the case (2), Lemma~\ref{prop:A.2} gives along the same lines of
argument that, almost surely,
\begin{eqnarray*}
\lim_{n\to\infty}\frac{1}n\sum
_{t=\hat k_{j-1}+1}^{\hat k_j}\rho_\tau(\hat
\varepsilon_t)
&>& \frac{1}{\lambda_j^*-\lambda_{j-1}^*} \Biggl\{\bigl(\lambda
_{j^\prime
}^0-
\lambda_{j-1}^*\bigr)\rho_{\tau,j^\prime}\\
&&\hphantom{\frac{1}{\lambda_j^*-\lambda_{j-1}^*} \Biggl\{}{}+\sum
_{\ell=j^\prime
+1}^{j^\prime+J+1}\bigl(\lambda_\ell^0-
\lambda_{\ell-1}^0\bigr)\rho_{\tau
,\ell}+\bigl(
\lambda_j^*-\lambda_{j^\prime+J}^0\bigr)
\rho_{\tau,j^\prime
+J+1} \Biggr\}.
\end{eqnarray*}
Taken together, these two inequalities, combined with the fact that
asymptotically all penalty terms in the definition of the \textsc{mdl}
in \eqref{mdl-lim} vanish, give, almost surely,
\begin{eqnarray*}
\lim_{n\to\infty}\frac{1}n{\textsc{mdl}}
\bigl(m^0,\hat\lambda,\hat p|\tau\bigr) &=&\lim_{n\to\infty}
\frac{1}n\sum_{j=1}^{m^0+1}\sum
_{t=\hat
k_{j-1}+1}^{\hat k_j}\rho_\tau(\hat
\varepsilon_t)
\\
&>&\lim_{n\to\infty}\frac{1}n\sum
_{j=1}^{m^0+1}\sum_{t=k^0_{j-1}+1}^{k_j^0}
\rho_\tau(\varepsilon_t) =\lim_{n\to\infty}
\textsc{mdl}\bigl(m^0,\lambda^0,p^0|\tau
\bigr),
\end{eqnarray*}
which is a contradiction to the definition of the MDL minimizer.
\end{pf*}
\begin{pf*}{Proof of Corollary~\ref{cor:1}}
Recall that the minimum description length criterion for multiple
quantiles $(\tau_{1},\ldots,\tau_{L})$ is given in (\ref
{eqn:qar11}). If follows from Theorem~\ref{th:1} that at any
individual quantile $\tau_{\ell}$, the minimizer, say, $(\hat
{\lambda}_\ell,\hat{p}_\ell)$ of the minimum description length
criterion (\ref{eqn:qar10}) is consistent for $(\lambda^0,p^0)$. It
follows that the minimizer $(\hat{\lambda},\hat{p})$ of (\ref
{eqn:qar11}) is consistent as it is a weighted sum of several criteria
in the form of (\ref{eqn:qar10}).
\end{pf*}

\end{appendix}



%

\printhistory
\end{document}